\documentclass[12pt]{amsart}
\usepackage{amsmath} 
\usepackage{latexsym} 
\usepackage{amssymb} 
\usepackage{graphicx}
\newtheorem{theorem}{Theorem} 
\newtheorem{proposition}{Proposition} 
\newtheorem{lemma}{Lemma}

\newtheorem{example}{Example} 

\newtheorem{remark}{Remark} 

\makeatletter
\@namedef{subjclassname@2020}{%
  \textup{2020} Mathematics Subject Classification}
\makeatother
\begin{document} 
%
\def\R {{\mathbb{R}}} 
\def\N {{\mathbb{N}}} 
\def\C {{\mathbb{C}}} 
\def\Z {{\mathbb{Z}}} 
\def\Q {{\mathbb{Q}}} 
\def\phi{\varphi} 
\def\epsilon{\varepsilon} 
\def\ma{{\mathcal A}} 
%
\def\tb#1{\|\kern -1.2pt | #1 \|\kern -1.2pt |}  
\def\Qed{\qed\par\medskip\noindent} 
%
\title[An interpolation problem in the Denjoy-Carleman classes ]{An interpolation problem in the Denjoy-Carleman classes}
\author{Paolo Albano}  
\address 
{Dipartimento di Matematica,  
Universit\`a di Bologna, Piazza 
di Porta San Donato 5, 40127 Bologna, Italy}  
\email{paolo.albano@unibo.it} 
\author{Marco Mughetti}  
\address 
{Dipartimento di Matematica,  
Universit\`a di Bologna, Piazza 
di Porta San Donato 5, 40127 Bologna, Italy}  
\email{marco.mughetti@unibo.it} 
   
\date{\today} 
 
\begin{abstract} Inspired by some iterative algorithms useful for proving the real analyticity (or the Gevrey regularity) of a solution of a linear partial differential equation with real-analytic coefficients, we consider the following question. Given a smooth function defined on $[a,b]\subset \R$ and given an increasing divergent sequence $d_n$ of positive integers such that the derivative of order $d_n$ of $f$ has a growth of the type $M_{d_n}$, when can we deduce that $f$ is a function in the Denjoy-Carleman class $C^M([a,b])$? We provide a positive result, and we show that a suitable condition on the gaps between the terms of the sequence $d_n$ is needed. 
\end{abstract} 
 
\subjclass[2020]{41A17, 26E05, 26E10}
\keywords{smooth functions, Denjoy-Carleman classes, quasi-analytic functions, Gevrey functions, analytic functions, interpolation} 
 
\maketitle

\section{Introduction and statement of the results} 
\setcounter{equation}{0} 
\setcounter{theorem}{0} 
\setcounter{proposition}{0}   
\setcounter{lemma}{0} 
\setcounter{corollary}{0}  
\setcounter{definition}{0} 

A way to prove the local regularity of the solutions of linear elliptic partial differential equations with real-analytic coefficients consists in the so called $L^2$--methods, i.e. an iterative procedure based on the use of the elliptic estimate (see e.g.  \cite[p.207]{BJS}). 
This approach extends to degenerate elliptic equations, such as sums of squares of vector fields, satisfying a subelliptic estimate (see, e.g., \cite{BM1,BM2} and \cite{BM3} for recent applications of this method). 
Then, a natural question is which derivatives one should control in order to conclude that a function is in a given Gevrey class (or, in particular, is real-analytic).  
This is the main motivation for the present paper. 
In order to put in evidence the essential points in our proofs, we state our results for the Denjoy-Carleman classes, and, for the sake of the simplicity, we limit our considerations to functions of one variable.  
To be definite:  let $M_0=1$, $M_1, \ldots $ be a sequence of positive numbers and consider
the Denjoy-Carleman class $C^M$
\begin{multline*}
C^M([a,b])=\{ f\in C^\infty ([a,b])\ | \ \exists K>0\text{ such that }
\\
|f^{(n)}(x)|\le K^{n+1} M_n\text{ if }x\in [a,b],\, n=0,1,\ldots   \}. 
\end{multline*}
In what follows we assume that the sequence $M_n$ is log-convex. That is, for  $j<\ell <k$, we have that  
\begin{equation}\label{eq:logcon} 
M_\ell \le M_j^{\frac{k-\ell }{k-j}} M_k^{\frac{\ell -j}{k-j}}.
\end{equation}
We denote by $\N$ the set of all the non-negative integers while $\N_+$ stands for the set of the positive integers. 
Condition \eqref{eq:logcon} can be equivalently stated as
\begin{equation}\label{eq:logcon2} 
M_n^2\le M_{n-1}M_{n+1},\qquad \forall n\in \N_+ , 
\end{equation}
i.e. $M_{n+1}/M_n$ is an increasing sequence. In particular,  if $M_n$ is log-convex then 
\begin{equation}\label{eq:monotonicity} 
M_n^{\frac 1n}\quad \text{ is an increasing sequence for }n\in \N_+ .   
\end{equation}
For the reader convenience, we provide the proofs  of the the assertions   
$\eqref{eq:logcon}\iff \eqref{eq:logcon2}$ and $\eqref{eq:logcon2} \implies \eqref{eq:monotonicity}$ in the Appendix \ref{A}.

We will assume the additional condition that there exists $m_0\geq 0$ such that  
\begin{equation}\label{eq:inter}
M_{j}\le M_k^{\frac jk}M_i^{\frac ji},\text{ for $i,k>m_0$ with $i<j$ and $j/i<k$.} 
\end{equation} 
\begin{example} 
(i) Our first example of a sequence satisfying \eqref{eq:logcon} and \eqref{eq:inter} is $M_n= n^{ns}$, $n\in \N_+ $, for a suitable $s\geq 1$.  
This choice corresponds to the Gevrey class $G^s$   (the set of all the real-analytic functions in the case $s=1$). We point out that \eqref{eq:inter} holds with $m_0=1$. 
  
\noindent   
(ii) Our second example is $M_0=M_1=1$, $M_2=\sqrt{M_3}$, $M_n= n^{ns_1}(\log n)^{ns_2}$ $(n=3,4,\ldots )$, for suitable $s_1,s_2\geq 1$. This sequence is log-convex and
satisfies \eqref{eq:inter} with $m_0=e^2$.
By the Denjoy-Carleman Theorem, $C^M([a,b])$ is a quasi-analytic class for $s_1=s_2=1$. (We recall that a class $C^M([a,b])$ is called quasi-analytic if for every $u\in C^M$ vanishing of infinite order at a point in $[a,b]$, it follows that $u$ is identically zero on $[a,b]$.)

We observe that, for $s_1>1$ or $s_2>1$, $C^M([a,b])$ is a non-quasi-analytic class (other than Gevrey). 

For the proof of the claims  in (i) and (ii), we refer the interested reader to the Appendix \ref{app2}. 
\end{example} 
Furthermore, we assume that the class $C^M$ contains the real-analytic functions on $[a,b]$, i.e. that there exists $c>0$ such that 
\begin{equation}\label{eq:anainc}
M_n\geq c^n n^n.  
\end{equation}

For a smooth function $f:[a,b]\longrightarrow \R$, we define 
\begin{equation}\label{fn}
F_n=\max_{x\in [a,b]}| f^{(n)}(x)|,\qquad n\in \N.
\end{equation}
Let $d_n>0$ be an increasing divergent sequence of integers such that, for a suitable $K>0$,  
\begin{equation}\label{eq:mest} 
F_{d_n}\le K^{d_n+1}M_{d_n}, \qquad \text{ for every }n\in \N.  
\end{equation}
 We consider the following problem: given a function $f\in C^\infty([a,b])$ satisfying \eqref{eq:mest},  under what condition on $\{d_n\}_n$ can we conclude that $f$ belongs to the class $C^M$? 

First, we observe that Condition \eqref{eq:mest} yields some additional properties on the function $f$ when $C^M$ is the class of real-analytic functions, i.e. $M_k=k^k$, $k\in \N_+$.  Precisely, we have the following 
\begin{proposition}\label{p}
Let $f\in C^\infty ([a,b])$ and let $C^F([a,b])$ be the corresponding Denjoy-Carleman class defined by the sequence $\{F_n\}_n$ 
in \eqref{fn}. Let $d_n>0$ be an increasing divergent sequence of integers. 
Assume that $F_{d_n}\le K^{d_n+1} d_n^{d_n}$ for a suitable $K>0$ and for every $n\in \N$.
Then, $f$ is identically zero, whenever $f$ vanishes of infinite order at a point of $[a,b]$. 
\end{proposition}
In other words, the class $C^F([a,b])$ is quasi-analytic (independently on the sequence $d_n$).
\begin{remark}
In general, we can not deduce that the class $C^F([a,b])$, given by Proposition \ref{p}, is contained in the one of the real-analytic functions, this a consequence of Theorem \ref{t2} with $M_n=n^n$.
\end{remark}
If we make an additional assumption on the growth of  $d_{n+1}-d_n$, we obtain the following 
\begin{theorem}\label{t}
Let $f\in C^\infty ([a,b])$,  assume \eqref{eq:logcon}, \eqref{eq:inter}, \eqref{eq:anainc}, \eqref{eq:mest}  and that there exists $c_0>0$ such that 
\begin{equation}\label{eq:gapc}
d_{n+1}/d_n\le c_0,\quad\text{ for }\quad n\in \N.
\end{equation}
 Then, $f\in C^M([a,b])$.   
\end{theorem}
In particular, in order for $f\in C^\infty$ to be in the Gevrey class $G^s$, it is enough to control only a set of derivatives $f^{(d_n)}$ where the sequence $\{d_n\}_n$ satisfies \eqref{eq:gapc}.
\begin{remark}
(i) We point out that the above result applies to both quasi-analytic and non-quasi-analytic Denjoy-Carleman classes. 

\noindent  
(ii) Theorem \ref{t} is related to the Carleman problem (see \cite{CM} and \cite{G}). 
\end{remark}
Finally, we show that Theorem \ref{t} may fail in the absence of  a condition on the sequence $\{d_n\}_n$.  
\begin{theorem}\label{t2}
Let  $\{M_n\}_n$ be an arbitrary sequence of positive numbers such that 
\begin{equation}\label{eq:anaincw}
\limsup_{n\to +\infty } M_n^{\frac 1n}=+\infty . 
\end{equation}
Then, we can find an increasing, divergent, sequence of positive integers $\{d_n\}_n$, and a function  $f\in C^\infty ([a,b])$, such that  \eqref{eq:mest} holds true but  $f\notin C^M([a,b])$. 
\end{theorem}

\section{Proofs} 
\setcounter{equation}{0} 
\setcounter{theorem}{0} 
\setcounter{proposition}{0}   
\setcounter{lemma}{0} 
\setcounter{corollary}{0}  
\setcounter{definition}{0} 

 \subsection{Proof of Proposition \ref{p}} 
 The proof of the Proposition \ref{p} is an elementary computation based on the Taylor expansion (see e.g. \cite{C}). 
Indeed, let $x_0\in [a,b]$, then by the Taylor formula, we find that 
$$
|f(x)|\le  K^{d_n+1} \frac{{d_n}^{d_n}}{d_n!} |x-x_0|^{d_n}\leq K(Ke|x-x_0|)^{d_n},\qquad \forall x\in [a,b], 
$$
and taking the limit as $n\to\infty$ in the formula above, we deduce that $f(x)=0$, for every $x\in [a,b]$ with $|x-x_0|<1/(Ke)$. 
Then, the conclusion follows from finitely many interations. 
 \subsection{Proof of Theorem \ref{t}} 

We want to show that $f\in C^M([a,b])$  if $f\in C^\infty ([a,b])$ and $\{M_n\}_n$ is a log-convex sequence such that 
\begin{itemize} 
\item[$(A)$] $M_{j}\le M_k^{\frac jk}M_i^{\frac ji}$, for $i,k>m_0$ with $i<j$ and $j/i<k$,  for a suitable $m_0\geq0$;
\item[$(B)$]  $M_n\geq c^n n^n$,  for a suitable $c>0$; 
\item[$(C)$] $F_{d_n}\le K^{d_n+1}M_{d_n}$ (see \eqref{fn}), 
\end{itemize}
 where  $\{d_n\}_n$ is an increasing divergent sequence of natural numbers satisfying the gap condition,  
\begin{equation}\label{eq:bound}
\frac{d_{n+1}}{d_n}\le c_0 ,\qquad \forall n\in \N,
\end{equation} 
 for a suitable positive integer $c_0$.

 The idea of the proof consists in showing that the control of the derivatives of length $d_n$ yields, by interpolation, an estimate for the intermediate derivatives (provided that the rescaled gap $(d_{n+1}-d_n)/d_n$ is bounded, uniformly w.r.t. $n$).   

For this purpose, we recall an estimate due to Cartan and Gorny (see e.g. \cite{C1,C2} and \cite{G}). 
\begin{lemma}\label{l} 
Let $g$ be a function $m$--times differentiable on the closed interval $[a,b]$ and set
$G_\ell=\max_{x\in [a,b]}| g^{(\ell)}(x)|$, for every $\ell\in\N$. Then, for every $m\in\N$, with $m\geq 2$, and 
$k\in \{ 1,\ldots ,m-1\}$, one has  
\begin{equation}\label{eq:cago}
G_k\le 2\left (  \frac{e^2m}{k}\right )^k G_0^{1-\frac km} \left ( \max \left \{  m!\, G_0 \left (\frac{2}{ b-a }\right )^m,\, G_m\right \} \right )^{\frac km}. 
\end{equation}
\end{lemma}

\vspace{0,5cm}
Let $d_n<\ell <d_{n+1}$ and choose $g=f^{(d_n)}$, $m=d_{n+1}-d_n$ and $k=\ell-d_n$ in \eqref{eq:cago}. 
It turns out that  
 $$
 \frac 1p:=1-\frac km=\frac{d_{n+1}-\ell }{d_{n+1}-d_n}\, ,  \qquad \frac 1q:=\frac km=\frac{\ell -d_n}{d_{n+1}-d_n}\, , 
 $$
and 
 \begin{equation}\label{eq:interl} 
 \frac{d_n}p+\frac{d_{n+1}}q=\ell . 
 \end{equation}
Whence \eqref{eq:cago} reads as
\begin{multline}\label{eq:gorny}
F_\ell \le  2   \left ( \frac{e^2 (d_{n+1}-d_n) }{\ell -d_n}  \right )^{\ell -d_n} \times
\\
  \max \left\{ \frac{     (d_{n+1} -d_n)!^{\frac 1q}   }{((b-a)/2)^{\ell -d_n} }   F_{d_n}  , F_{d_n}^{\frac 1p} F_{d_{n+1}}^{\frac 1q}      \right\}. 
\end{multline} 
Now, in view of \eqref{eq:bound}, we have that 
\begin{multline*} 
\displaystyle{ \left ( \frac{e^2 (d_{n+1}-d_n) }{\ell -d_n}  \right )^{\ell -d_n}= \left ( e^2 \left (  1+\frac{  d_{n+1}-\ell }{\ell -d_n}\right )   \right )^{\ell -d_n}\leq e^{2(\ell-d_n)+d_{n+1}-\ell}}
\\
\displaystyle{\le e^{\ell +d_n \left ( \frac{d_{n+1}}{d_n} -2\right )}\le e^{\ell (c_0-1)}}.
 \end{multline*} 
Furthermore, 
$$ 
 \displaystyle{\frac{ (d_{n+1} -d_n)!^{\frac 1q}}{((b-a)/2)^{\ell -d_n} }\le\frac{ (d_{n+1} -d_n)^{\frac {d_{n+1}-d_n}q}}{((b-a)/2)^{\ell -d_n} }\le \left (\frac{2d_n( c_0-1) }{b-a}\right )^{\ell -d_n}}\, , 
$$
and, since $F_{d_n}\le K^{d_n+1}M_{d_n}$, we obtain that
\begin{equation}\label{1}
F_\ell \le   
 \max \left\{ 2e^{\ell (c_0-1)} \left (\frac{2d_n  ( c_0-1)    }{b-a}\right )^{\ell -d_n} K^{d_n+1}M_{d_n}  ,   2e^{\ell (c_0-1)} K^{\frac{d_n+1}p} M_{d_n}^{\frac 1p} K^{\frac{d_{n+1}+1}q} M_{d_{n+1}}^{\frac 1q}      \right\}. 
\end{equation}
We observe that, by \eqref{eq:interl} and the fact that $d_n<\ell $, we have that  
\begin{equation}\label{eq:aggn} 
K^{\frac{d_n+1}p} K^{\frac{d_{n+1}+1}q}=K^{\ell +1},\text{ and }K^{d_n+1}\le K^{\ell +1}
\end{equation} 
Moreover, by Assumption $(A)$ above\footnote{ We observe that, without loss of generality, we can take $c_0$ and $d_n$ larger than the number $m_0$ in \eqref{eq:inter}} and by the gap condition \eqref{eq:bound}, we see that
$$
M_{d_n}^{\frac 1p} M_{d_{n+1}}^{\frac 1q} \le M_{d_n}^{\frac 1p}M_{c_0}^{\frac{d_{n+1}}{c_0 q}} M_{d_n}^{ \frac{d_{n+1}}{qd_n}}.
$$
Whence, as a consequence of the log-convexity of $\{M_n\}_n$, \eqref{eq:monotonicity}, and \eqref{eq:interl}, we find that
\begin{equation}\label{2}
M_{d_n}^{\frac 1p} M_{d_{n+1}}^{\frac 1q}
\le M_{\ell}^{\frac{d_n}{\ell p}} M_{c_0}^{\frac{d_{n+1}}{c_0q}} M_{\ell }^{ \frac{d_{n+1}}{q\ell }}
\le M_{\ell} M_{c_0}^{\frac{\ell }{ c_0}}. 
\end{equation}
Assumption $(B)$ yields that 
\begin{equation}\label{3}
d_n^{\ell -d_n}\le \left (\frac{M_{d_n}^{\frac{ 1}{d_n}} }{c}\right )^{\ell -d_n} .
\end{equation}
Hence, taking together \eqref{1}, \eqref{eq:aggn}, \eqref{2}, \eqref{3} and  using the  log-convexity of $\{M_n\}_n$, \eqref{eq:monotonicity}, we deduce that 
\begin{multline*}
F_\ell \le  \max \left\{ 2e^{\ell (c_0-1)} \left (\frac{2M_{d_n}^{\frac{ 1}{d_n}} ( c_0-1)    }{c(b-a)}\right )^{\ell -d_n} M_{d_n}  ,\,   2e^{\ell (c_0-1)}        M_{\ell} M_{c_0}^{\frac{\ell}{c_0}}\right\} K^{\ell +1}
\\
\le  \max \left\{ 2e^{\ell (c_0-1)} \left (\frac{2(c_0-1)}{c(b-a))}\right )^{\ell -d_n} , \, 2e^{\ell (c_0-1)}  M_{c_0}^{\frac{\ell }{c_0}}\right\} M_{\ell} 
\le C_1^\ell K^{\ell +1} M_{\ell},   
\end{multline*}
 where 
 $$
 C_1=  2e^{c_0-1} \left (\frac{2(c_0-1)}{c(b-a)}+ M_{c_0}^{\frac 1{c_0}}\right ).   
 $$
 Then,  we conclude that there exists $K_1>0$ such that   
 $$
 F_n\le K_1^{n+1}M_n ,\quad \text{  for every }n\in \N ,
 $$
 i.e.  $f\in C^M([a,b])$. 
 This completes our proof of Theorem \ref{t}.

\subsection{Proof of Theorem \ref{t2}} 

The proof is based on the following result, which ensures the existence, in every class $C^N$, of a function that attains the bounds $N_j$.
Although, different formulations are already present in the literature (see, e.g., \cite{T} for the case of complex valued functions), for the sake of completeness we provide its proof.    
\begin{lemma}\label{l2}
Let $\{N_j\}_j$ be a positive sequence satisfying \eqref{eq:logcon}. Then, there exists $f\in C^N([a,b])$ such that 
\begin{equation}\label{eq:optimal}
\left |f^{(j)}\left (\frac{a+b}2\right )\right |\geq N_j,\qquad \forall j\in \N. 
\end{equation} 
\end{lemma} 
 \begin{proof}
 Set $m_j=N_{j+1}/N_j$ and observe that, by \eqref{eq:logcon}, it is an increasing sequence. 
 We note that 
 \begin{equation}\label{eq:partial}
 \left ( \frac{1}{m_k}  \right )^{k-j} \le \frac{N_j}{N_k}\qquad \forall k,j\in \N. 
 \end{equation}
Indeed, \eqref{eq:partial} trivially holds in the case of $j=k$. Moreover, if $j<k$, we have that 
$$
\left ( \frac{1}{m_k}  \right )^{k-j} \le \frac{1}{m_{k-1}\cdots m_{k-(k-j)}}=\frac{N_j}{N_k}\, .
$$  
Finally, if $j>k$, we obtain 
$$
m_k^{j-k} \le m_{k}\cdots m_{j-1}=\frac{N_j}{N_k}\, . 
$$ 
Summing up,  \eqref{eq:partial} holds true.

Let us define 
$$
g(x)=\sum_{k=0}^{\infty } \frac{N_k}{(2m_k)^k} (\cos (2m_kx)+\sin (2m_k x)). 
$$
We notice that $g\in C^N(I)$ for every interval $I\subset \R$.  Indeed, in light of  \eqref{eq:partial},  we get 
$$
|g^{(n)}(x)|\le \sum_{k=0}^{\infty } \frac{N_k}{(2m_k)^k} (2m_k)^n 2\le   N_n\sum_{k=0}^{\infty } \frac{2^{n+1}}{2^k} \le 2^{n+2}N_n.  
$$
Furthermore, we have that 
$$
|g^{(n)}(0)|=\left | \left( \sum_{k=0}^\infty  \frac{N_k}{(2m_k)^{k-n}} \right) (\cos^{(n)}(0)+ \sin^{(n)}(0)) \right |,
$$
and, since 
$$
\begin{cases}
\cos^{(n)}(0)=(-1)^{n/2}\text{ and } \sin^{(n)}(0)=0,\quad \text{ for $n$ even,} 
\\
\cos^{(n)}(0)=0\text{ and } \sin^{(n)}(0)=(-1)^{\frac{n-1}{2}},\quad \text{ for $n$ odd,}
\end{cases} 
$$
we deduce that 
$$
|g^{(n)}(0)|=\sum_{k=0}^\infty  \frac{N_k}{(2m_k)^{k-n}}    \geq N_n
$$
The conclusion follows by taking 
$$
f(x):=g\left ( x-\frac{a+b}2 \right ),\qquad (x\in [a,b]).  
$$
\end{proof} 
Then, the proof of Theorem \ref{t2} reduces to show that given a positive sequence $\{M_n\}_n$ satisfying \eqref{eq:anaincw},  there exist 
\begin{itemize}
\item a sequence $\{N_n\}_n$ satisfying  \eqref{eq:logcon},
\item two divergent sequences  of positive integers $\{d_n\}_n$ and $\{i_n\}_n$, 
\end{itemize}
so that   
$$
N_{d_n}= M_{d_n},\qquad n=0,1,\ldots  
$$
and
\begin{equation}\label{eq:fg}
 \frac{N_{i_n}}{M_{i_n}}=2^{2^{i_n}} . 
\end{equation}    
We point out that, once the existence of such a sequence $\{N_n\}_n$ is established, Lemma \ref{l2} yields the existence of a function $f\in C^N([a,b])$ such that $N_{d_n}\le M_{d_n}$. On the other hand,  due to  \eqref{eq:fg} , we have that $f\notin C^M([a,b])$.   
In other words, in general (if the gap 
$(d_{n+1}-d_n)/d_n$ is suitably large), the interpolation between $d_n$ and $d_{n+1}$ does not provide the bounds on the derivatives of $f$ ensuring that $f\in C^M([a,b])$.

We construct inductively the sequences $\{N_n\}_n$, $\{d_n\}_n$ and $\{i_n\}_n$. 
For this purpose, it is useful to define 
$$
m_n=\frac{N_{n+1}}{N_n}. 
$$
We observe that $N_n$ satisfies \eqref{eq:logcon} if and only if $m_n$ is an increasing sequence; furthermore, without loss of generality, we may assume that  $M_0=N_0=1$.  
Clearly, we have
$$
m_0\cdot...\cdot m_{j-1}=\frac{N_j}{N_0}=N_j.
$$
Let us fix an arbitrary positive integer $i_0$. 
We claim that there exists  a positive integer $d_0>i_0$  so that setting 
\begin{equation}\label{eq:step1}
\begin{cases}
m_0=\ldots =m_{i_0-1}=(2^{2^{i_0}}M_{i_0})^{\frac{1}{i_0}},
\\ 
m_{i_0}=\ldots =m_{d_0-1}=\left (\frac{M_{d_0}}{2^{2^{i_0}}M_{i_0}}\right )^{\frac 1{d_0-i_0}}
\end{cases}
\end{equation}
we have
\begin{equation}\label{eq:monstep1}
m_j\le m_{j+1},\qquad \text{ for }j=0,\ldots ,d_0-2 . 
\end{equation}
We observe that \eqref{eq:monstep1} is equivalent to the inequality
$$
(2^{2^{i_0}}M_{i_0})^{\frac{1}{i_0}}\le \left (\frac{M_{d_0}}{2^{2^{i_0}}M_{i_0}}\right )^{\frac 1{d_0-i_0}}\iff
(2^{2^{i_0}}M_{i_0})^{\frac{1}{i_0}} \le M_{d_0}^{\frac{1}{d_0}}.
$$
In view of \eqref{eq:anaincw},  for every $c>0$ there exists $d_0$ such that  $M_{d_0}^{\frac{1}{d_0}}\geq c$. Hence, we find that 
$$
M_{d_0}^{\frac{1}{d_0}}\geq c\geq (2^{2^{i_0}}M_{i_0})^{\frac{1}{i_0}},
$$
provided  that $c>  (2^{2^{i_0}}M_{i_0})^{\frac{1}{i_0}}$, i.e.  \eqref{eq:monstep1} holds true.\\ 
We point out that, as a consequence of our construction, we have
\begin{equation}\label{eq:m1}
\begin{cases}
N_j=(2^{2^{i_0}}M_{i_0})^{\frac j{i_0}},\qquad &\text{for } j=0,\ldots ,i_0, 
\\
 N_j= 2^{2^{i_0}}M_{i_0}  \left (\frac{M_{d_0}}{2^{2^{i_0}}M_{i_0}}\right )^{\frac {j-i_0}{d_0-i_0}},\,  &\text{for } j=i_0+1,\ldots d_0, 
\end{cases} 
\end{equation}
\begin{equation}\label{eq:mm1}
\frac{N_{i_0}}{M_{i_0}}=2^{2^{i_0}},\quad\textrm{and}\quad N_{d_0}=M_{d_0}.
\end{equation}
 Let us show that there exist two positive integers $d_1>i_1>d_0$ so that setting 
\begin{equation}\label{eq:step2}
\begin{cases}
m_{d_0}=\ldots =m_{i_1-1}=\left (\frac{2^{2^{i_1}}M_{i_1}}{M_{d_0}}\right )^{\frac{1}{i_1-d_0}},
\\ 
m_{i_1}=\ldots =m_{d_1-1}=\left (\frac{M_{d_1}}{2^{2^{i_1}}M_{i_1}}\right )^{\frac{1}{d_1-i_1}},
\end{cases}
\end{equation}
we have that 
\begin{equation}\label{eq:monstep2}
m_j\le m_{j+1},\qquad \text{ for }j=0,\ldots ,d_1-2 . 
\end{equation}
 We observe that \eqref{eq:monstep2} is equivalent to 
 \begin{equation}\label{eq:monstep2bis}
 \begin{cases}
 \left (\frac{2^{2^{i_1}}M_{i_1}}{M_{d_0}}\right )^{\frac{1}{i_1-d_0}}\geq \left (\frac{M_{d_0}}{2^{2^{i_0}}M_{i_0}}\right )^{\frac 1{d_0-i_0}},
 \\
 \left (\frac{M_{d_1}}{2^{2^{i_1}}M_{i_1}}\right )^{\frac{1}{d_1-i_1}}\geq \left (\frac{2^{2^{i_1}}M_{i_1}}{M_{d_0}}\right )^{\frac{1}{i_1-d_0}}. 
 \end{cases}
 \end{equation}
 As for the first inequality in \eqref{eq:monstep2bis}, it is enough to show that 
 $$
 \left (\frac{2^{2^{i_1}}M_{i_1}}{M_{d_0}}\right )^{\frac{1}{i_1-d_0}}\geq \left (M_{d_0}\right )^{\frac 1{d_0-i_0}}\, ,
 $$
which is equivalent to 
 $$
 2^{2^{i_1}}M_{i_1}\geq M_{d_0}^{\frac{i_1-i_0}{d_0-i_0}}.
 $$
Once more, by \eqref{eq:anaincw}, this last inequality can be satisfied provided  $i_1>d_0$ is chosen suitably large. 
 Furthermore, the second inequality in \eqref{eq:monstep2bis} can be rewritten as 
 $$
 M_{d_1}\geq \frac{(2^{2^{i_1}}M_{i_1})^{\frac{d_1-d_0}{i_1-d_0}} }{M_{d_0}^{\frac{d_1-i_1}{i_1-d_0}}}\, , 
 $$
 which can be fulfilled, due to \eqref{eq:anaincw}, provided that $d_1>i_1$ be large enough. \\
Summing up, we have constructed a sequence $N_j$, $j=0,\ldots d_1$, such that $m_j=\frac{N_{j+1}}{N_j}$ is an increasing sequence for $j=0,...,d_1$, 
 and
$$
N_{d_0}=M_{d_0},\quad N_{d_1}=M_{d_1},\quad N_{i_0}/M_{i_0}=2^{2^{i_0}}\quad \textrm{and}\quad N_{i_1}/M_{i_1}=2^{2^{i_1}},
$$ 
for suitable positive integers $0<i_0<d_0<i_1<d_1$. \\
Now, let us suppose that we have already defined $d_0<i_0<\cdots <i_{n}< d_n$, and 
$N_j$, $j=1,\ldots d_{n}$, such that 
$$
\begin{cases}
m_j=\frac{N_{j+1}}{N_j}\quad \textrm{are increasing for } j=0,...,d_n-1,
\\
N_{d_k}=M_{d_k}, \qquad k=0,\ldots ,n,
\\
\frac{N_{i_k}}{M_{i_k}}=2^{2^{i_k}},\qquad k=0,\ldots ,n. 
\end{cases} 
$$
In order to complete the proof of Theorem \ref{t2} it suffices to show that we can find $i_{n+1}$, $d_{n+1}$ and $ N_j$, $j=d_{n}+1,\ldots,d_{n+1}$   such that 
\begin{enumerate}
\item $d_n<i_{n+1}<d_{n+1}$, 
\item $m_j=\frac{N_{j+1}}{N_j}$ are increasing for $j=0,\ldots , d_{n+1}$, 
\item $N_{d_{n+1}}=M_{d_{n+1}}$ and $\frac{N_{i_{n+1}}}{M_{i_{n+1}}}=2^{2^{i_{n+1}}}$. 
\end{enumerate} 
Set 
\begin{equation}\label{eq:stepn}
\begin{cases}
m_{d_n}=\ldots =m_{i_{n+1}-1}=\left (\frac{2^{2^{i_{n+1}}} M_{i_{n+1}}}{M_{d_n}}\right )^{\frac{1}{i_{n+1}-d_n}},
\\ 
m_{i_{n+1}}=\ldots =m_{d_{n+1}-1}=\left (\frac{M_{d_{n+1}}}{2^{2^{i_{n+1}}} M_{i_{n+1}}}\right )^{\frac{1}{d_{n+1}-i_{n+1}}},
\end{cases}
\end{equation}
we have that 
\begin{equation}\label{eq:monstepn}
m_j\le m_{j+1},\qquad \text{ for }j=0,\ldots ,d_{n+1}-2 . 
\end{equation}
We observe that \eqref{eq:monstepn} reduces to verifying 
\begin{equation}\label{eq:monstepnbis}
 \begin{cases}
 \left (\frac{2^{2^{i_{n+1}}} M_{i_{n+1}}}{M_{d_n}}\right )^{\frac{1}{i_{n+1}-d_n}}\geq m_{d_n-1} ,
 \\
   \left (\frac{M_{d_{n+1}}}{2^{2^{i_{n+1}}} M_{i_{n+1}}}\right )^{\frac{1}{d_{n+1}-i_{n+1}}}\geq \left (\frac{2^{2^{i_{n+1}}}M_{i_{n+1}}}{M_{d_n}}\right )^{\frac{1}{i_{n+1}-d_n}}. 
 \end{cases}
 \end{equation}
We observe that the first inequality in \eqref{eq:monstepnbis} is equivalent to 
$$
  M_{i_{n+1}}\geq \frac{M_{d_n} m_{d_n-1}^{i_{n+1}-d_n}}{ 2^{2^{i_{n+1}}}} 
$$
which, by \eqref{eq:anaincw}, can be satisfied provided that $i_{n+1}>d_n$ be large enough. 
Finally, the second inequality in \eqref{eq:monstepnbis}, can be rewritten as 
$$
 M_{d_{n+1}} \geq \frac{( 2^{2^{i_{n+1}}} M_{i_{n+1}})^{\frac{d_{n+1}-d_{n}}{i_{n+1}-d_n}}    }{M_{d_n}^{\frac{d_{n+1}-i_{n+1}}{i_{n+1}-d_n}}    } 
$$
which, once more by \eqref{eq:anaincw}, can be satisfied provided that $d_{n+1}>i_{n+1}$ be large enough. 
Finally, as a consequence of \eqref{eq:stepn} we have:
$$
\frac{N_{i_{n+1}}}{N_{d_n}}=m_{d_n}\cdot m_{d_n+1}\cdot ... \cdot m_{i_{n+1}-1}=\frac{2^{2^{i_{n+1}}} M_{i_{n+1}}}{M_{d_n}}
$$
$$
\frac{N_{d_{n+1}}}{N_{d_n}}=m_{d_n}\cdot m_{d_n+1}\cdot ... \cdot m_{d_{n+1}-1}=\frac{M_{d_{n+1}}}{M_{d_n}}
$$
Since $N_{d_n}=M_{d_n}$, this completes our proof of Theorem \ref{t2}.

\appendix
\section{} \label{app1}
 \setcounter{equation}{0} 
\setcounter{theorem}{0} 
\setcounter{proposition}{0}   
\setcounter{lemma}{0} 
\setcounter{corollary}{0}  
\setcounter{definition}{0}

Let $M_0=1$, $M_1,M_2,\ldots $ be a sequence of positive numbers. 
In this section we consider the following assertions  
\begin{itemize} 
\item[(A)] For $j<\ell <k$, $M_\ell \le M_j^{\frac{k-\ell }{k-j}} M_k^{\frac{\ell -j}{k-j}}$;
\item[(B)] $M_n^2\le M_{n-1}M_{n+1}$, for $n\in \N_+$; 
\item[(C)] $M_n^{\frac 1n}$ is an increasing sequence for $n\in \N_+$.   
\end{itemize} 
We show that $(A)\iff (B)$ and $(B)\implies (C)$.  

We observe that (B) can be rephrased by requiring that $M_{n}/M_{n-1}$ is an increasing sequence for $n\in \N_+$. 
Indeed, dividing both the sides of the inequality in (B) by $M_{n-1}M_{n}$,  we deduce that (B) holds if and only if $M_{n}/M_{n-1}$ is an increasing sequence for $n\in \N_+$. 

Now, let us show that $(B)\implies (C)$. 
(B) implies that 
$$
\frac{M_{n+1}}{M_{n}}\geq \frac{M_{n}}{M_{n-1}}\geq \frac{M_{n-1}}{M_{n-2}}\geq \ldots \geq\frac {M_1}{M_0}=M_1. 
$$
In particular, since $M_{n+1}/M_{n}$ is greater or equal than each of the factors $M_{n}/M_{n-1}$, $M_{n-1}/M_{n-2},\ldots ,M_1$,  we find that 
it is greater or equal than the geometrical mean of these factors, i.e.  
$$
 \frac{M_{n+1}}{M_{n}}\geq \left ( \frac{M_{n}}{M_{n-1}}\,  \frac{M_{n-1}}{M_{n-2}}\ldots M_1\right )^{\frac 1n}=M_n^{\frac 1n} . 
$$
Then, we find that $M_{n+1}\geq M_n^{\frac{n+1}{n}}$, i.e.  $M_n^{\frac 1n}$ is an increasing sequence. This completes the proof of $(B)\implies (C)$. 

Now, let consider the implication $(A)\implies (B)$.

Choosing $j=n-1$, $\ell =n$ and $k=n+1$, the inequality in (A) can be rewritten as $M_n \le M_{n-1}^{\frac 12} M_{n+1}^{\frac 12}$. 
Then, taking the square of both sides of the above inequality we deduce that $(A)\implies (B)$.

In order to show that $(B)\implies (A)$ we need the following 

\begin{lemma}\label{l:a} 
For every $n\in \N$ the sequence 
$$
\N_+\ni h\mapsto \left ( \frac{M_{n+h}}{M_n}\right )^{\frac 1h}\text{ is increasing.} 
$$
\end{lemma}
\begin{remark} 
The geometrical content of the above result is that a function is convex if and only if its slope is increasing. 
\end{remark} 
\begin{proof} 
The proof of the lemma proceeds by induction on $h$. We observe that 
$$
\left ( \frac{M_{n+h+1}}{M_n}\right )^{\frac 1{h+1}}\geq \left ( \frac{M_{n+h}}{M_n}\right )^{\frac 1h}
$$
is equivalent to the inequality 
\begin{equation}\label{eq:mpo}
M_{n+h}\le  M_{n}^{\frac 1{h+1}}M_{n+h+1}^{\frac h{h+1}}  ,  
\end{equation}
for every $h\in \N_+$.

For $h=1$, \eqref{eq:mpo} reduces to (B). 
Now, let us suppose that the inequality \eqref{eq:mpo} is satisfied for a suitable natural number $h$ and let us show that it is satisfied for $h+1$. 
Indeed, (B) implies that 
\begin{equation}\label{eq:ind1}
M_{n+h+1}\le M_{n+h}^{\frac 12 }M_{n+h+2}^{\frac 12 }.
\end{equation}
 Furthermore, by the inductive assumption, we have that 
 $$
 M_{n+h}^{\frac 12}\le M_{n}^{\frac 1{2(h+1)}}M_{n+h+1}^{\frac h{2(h+1)}}  
 $$
 and plugging the above inequality in \eqref{eq:ind1}, we find 
 $$
 M_{n+h+1}\le M_{n}^{\frac 1{2(h+1)}}M_{n+h+1}^{\frac h{2(h+1)}}  M_{n+h+2}^{\frac 12 }.  
 $$
 This last inequality yields that 
 $$
 M_{n+h+1}^{\frac {h+2}{2(h+1)}}\le M_{n}^{\frac 1{2(h+1)}}M_{n+h+2}^{\frac 12 }, 
 $$
 i.e. 
 $$
 M_{n+h+1}\le  M_{n}^{\frac 1{h+2}}M_{n+h+2}^{\frac {h+1}{h+2}}  
 $$
 and the proof of the lemma is completed. 
\end{proof}   
We observe that (A) can be rewritten as 
 \begin{itemize} 
 \item[(A$'$)] for $j<\ell <k$, $\log M_\ell  \le \frac{k-\ell}{k-j}  \log M_j +\frac{\ell -j}{k-j} M_k$.
 \end{itemize}
Now, Lemma \ref{l:a} implies that 
for every $n\in \N$ the sequence 
$$
\N_+\ni h \mapsto \frac{\log (M_{n+h})-\log (M_n)}{h}\quad \text{ is increasing.} 
$$
Then, for $j<\ell <k$, we find that 
$$
 \frac{\log (M_{\ell })-\log (M_j)}{\ell -j}\le  \frac{\log (M_{k})-\log (M_j)}{k-j} ,
$$
 i.e. $(A')$ holds. This completes the proof of the equivalence $(A)\iff (B)$.

 \section{} \label{app2}
\setcounter{equation}{0} 
\setcounter{theorem}{0} 
\setcounter{proposition}{0}   
\setcounter{lemma}{0} 
\setcounter{corollary}{0}  
\setcounter{definition}{0}

In this section we provide all the computations needed to justify the claims done in Example 1.

\noindent (i) 
Let $s\geq 1$ and $M_n=n^{ns}$, for $n\in \N_+$. 
In order to show that $M_n$ is log-convex it suffices to show that the function $f(x)=xs\log x$ is convex for $x\geq 1$ (it is clear that, due to the limitation $s\geq 1$, the factor $s$ here is immaterial). Since $f^{(2)}(x)=s/x>0$, for $x>0$, we deduce that $M_n$ is log-convex. 

Let us verifty that  $M_n$ satisfies \eqref{eq:inter}, i.e.  that there exists $m_0\geq 0$ such that  
$$
M_{j}\le M_k^{\frac jk}M_i^{\frac ji},\text{ for $i,k>m_0$ with $i<j$ and $j/i<k$.} 
$$ 
In this case, \eqref{eq:inter} can be written as 
$$
j^{js}\le k^{js} i^{js},\text{ for $i,k>m_0$ with $i<j$ and $j/i<k$,} 
$$
i.e. 
$$
j\le k i,\text{ for $i,k>m_0$ with $i<j$ and $j/i<k$,} 
$$
which trivially holds true with $m_0=1$.

\noindent (ii) 
Let $s_1,s_2\geq 1$, $M_0=M_1=1$, $M_2=\sqrt{M_3}$, and $M_n=n^{ns_1} (\log n)^{ns_2}$ $(n=3,4,\ldots )$. We claim that $M_n$ is log-convex. 
Recalling that the log-convexity of $M_n$ is equivalent to the fact that the sequence $M_{n+1}/M_n$ is increasing, in order to prove that $M_n$ is log-convex it suffices to verify that 
$M_0'=M_1'=1$, $M_2'=3^{3s_1/2}$, $M_n'=n^{n s_1}$  $(n=3,4,\ldots )$ and $M_0''=M_1''=1$, $M_2''=(\log 3)^{3s_2/2}$, $M_n''=(\log n)^{n s_2}$  $(n=3,4,\ldots )$ are separately log-convex. In light of (i) above, $M_n'$ is log-convex and the proof reduces to show that $M_n''$ has the same regularity. 
We have that 
$$
\frac{M_1''}{M_0''}\le \frac{M_2''}{M_1''}\iff 1\le (\log 3)^{3s_2}, 
$$
$$
\frac{M_2''}{M_1''}\le \frac{M_3''}{M_2''}\iff  M_3''\le M_3'', 
$$
$$
 \frac{M_3''}{M_2''}\le  \frac{M_4''}{M_3''}\iff M_3''^3=(\log 3)^9 \le M_4''^2=(\log 4)^8. 
$$
In  order to prove this last inequality it suffices to show that 
\begin{equation}\label{eq:log}
\log 3\le \log 4 /\log 3.
\end{equation}  
Using the concavity of the logarithm, we have that 
\begin{equation}\label{eq:logconc}
\log 3\le \log e +\frac 1e (3-e)=\frac 3e \le \frac 3{ 2+1/2}=\frac 65 ,  
\end{equation}  
and 
$$
\log 3\le \log 4 -\frac 14 \iff \frac{\log 4}{\log 3}\geq 1 +\frac{1}{4\log 3}. 
$$
The inequality above and \eqref{eq:logconc} yield 
$$
\frac{\log 4}{\log 3}\geq 1 +\frac{5}{24}\, \left (> 1+\frac 15 \geq \log 3  \right ),  
$$
and \eqref{eq:log} follows. 
It remains to show that $(\log n)^{ns_2}$ is log-convex for $n\geq 3$. Once more, this fact can be reduced to a one-variable problem: it suffices to verify that the function $f(x)=xs_2 \log (\log (x))$ is convex for $x\geq 3$. 
Now, $f^{(2)}(x)=s_2 (\log x -1)/[(\log x )^2x]>0$, for $x>e$, and the log-convexity of $M_n''$ (and $M_n$) follows.  

It remains to show that $M_n$ satisfies \eqref{eq:inter}. Clearly, it suffices to verify that $n^{ns_1}$ and $(\log n)^{ns_2}$ satisfy separately \eqref{eq:inter}. The factor  $n^{ns_1}$ was already treated in (i) above. Then, the proof reduces to show that 
$$
(\log j)^{js_2}\le (\log k)^{js_2}(\log i)^{js_2} ,\text{ for $i,k>m_0$ with $i<j$ and $j/i<k$,}
$$
i.e. 
\begin{equation}\label{eq:AB}
\log j\le \log k \log i ,\text{ for $i,k>m_0$ with $i<j$ and $j/i<k$.}
\end{equation} 
Now, the inequality $j/i<k$ and the monotonicity of the logarithm imply that 
$$
\log j < \log i + \log k . 
$$
Now, taking $m_0=e^2$, we have that 
$$
\log i + \log k \le \log i \log k , \text{ for $i,k>m_0$},  
$$
then \eqref{eq:AB}, and this completes the proof that $M_n$ satisfies \eqref{eq:inter}.

Finally, the fact that the class $C^M([a,b])$ is quasi-analytic for $s_1=s_2=1$ while it is non-quasi-analytic for $s_1>1$ or $s_2>1$,  follows by the Denjoy-Carleman Theorem (see e.g. \cite{C}). Indeed,  since $M_n^{\frac 1n}$ is increasing (in view of the log-convexity of $M_n$), the Denjoy-Carleman Theorem can be stated as follows:

the class $C^M([a,b])$ is quasi-analytic if and only if 
\begin{equation}\label{eq:DC}
\sum_{n=1}^\infty  \frac{1}{M_n^{\frac 1n}}=+\infty . 
\end{equation}
Hence, since $\sum_{n=1}^\infty  \frac{1}{M_n^{\frac 1n}}$ behaves as 
$$
\sum_{n=3}^\infty \frac1{n^{s_1} (\log n)^{s_2}}, \qquad (s_1,s_2\geq 1),  
$$
using the generalized integral associated with the series above, we conclude that the series in \eqref{eq:DC} diverges if and only if $s_1=s_2=1$. 
 
We point out that, by Lemma \ref{l2}, the class  $C^M([a,b])$ is a proper sub-class of $G^{s_1+\epsilon}$, for every $\epsilon >0$.

\section*{Declarations}
\begin{itemize}
\item{\bf Funding:} This work was partly supported by the National Group for Mathematical Analysis, Probability and Applications (GNAMPA) of the Italian Istituto Nazionale di Alta Matematica ``Francesco Severi''.
\end{itemize}

\end{document}